# Quasi-Distribution Appraisal Based on Piecewise Bézier Curves: An Objective Evaluation Method about Finite Element Analysis

Runkai. Wen[2,b]; Yukun. Chai[2,c]; Lingxin. Wang[3,d]; Ruochen. Du[2,e]; Xingtian. Long[4,f]; Zhiyang. Liu[3,g]; Peng. Wu[5,h]; Yiduo. Wang[1,a,*]

[1]Department of Mathematics, Beijing University of Chemical Technology, China.
[2]School of Internation Education, Beijing University of Chemical Technology, China.
[3]School of Civil and Hydraulic Engineering, Chongqing University of Science & Technology, China.
[4]School of Civil Engineering and Environment, Nanyang Technology University, Singapore
[5]School of Informatics, University of Edinburgh, UK.
[a]wangyd@mail.buct.edu.cn, [b]wenrunkai2023@gmail.com, [c]chaiyukun2024@163.com, [d]wanglinxin2001@163.com, [e]duruochen@qq.com, [f]longxingtian550@gmail.com, [g]liuzhiyang0519@163.com, [h]wpeng202409@163.com
*Corresponding author: Yiduo Wang.

**Abstract**
A class of quasi-distribution evaluation criteria based on piecewise Bézier curves is proposed to address the issue of the inability to objectively evaluate finite element models. During the optimization design of mechanical parts, finite element modeling is performed on their stress deformation, and the mesh node shape variable values are converted into distribution histogram data for piecewise Bézier curve fitting. Being dealt with area normalization method, the fitting curve could be regarded as a kind of probability density function (PDF), and its variance could be used to evaluate the finite element modeling results. The situation with the minimum variance is the optimal choice for overall deformation. Numerical experiments have indicated that the new method demonstrated the intrinsic characteristics of the finite element models of difference mechanical parts. As an objective appraisal method for evaluating finite element models, it is both effective and feasible.

## 1 Introduction

As is well known, the finite element method (FEM) is a highly efficient and commonly used computational method, which was developed based on the variational principle decades ago. With the widespread application of computer technology, the finite element method (FEM) has become an advanced CAE technology, widely used in various engineering fields to solve complex design and analysis problems, and has become an important tool in engineering design and analysis. This technique involves discretizing the domain under investigation into a series of smaller, easier-to-manage components using finite elements and subsequently reconnecting these elements at specific locations referred to as nodes[1]. FEM can also be used to calculate the internal stress, strain, and other parameters of an object. It can be combined with material properties and boundary conditions to analyze the static equilibrium state of the object structure[2], or to determine the circumferential mode frequencies of circular cylindrical shells of stator in electric motor[3]. As an engineering analysis method used to study the behavioral characteristics of contact or friction between objects, contact analysis often utilizes FEM to model the contact surfaces to describe the interaction between the contact surfaces. This approach solves many engineering problems, such as designing and evaluating the contact performance of mechanical components[4], assessing machining quality[5], estimating fatigue life[6], demonstrating cutting mechanism[7]. Increasingly, FEM has been broadly used in many fields, and this trend is dictated by the requirement of following: composite material analysis is an engineering analysis method primarily used to study the properties and behavioral characteristics of composite materials, which are composed of two or more materials. This method describes the physical and mechanical properties of materials by establishing mathematical models in order to predict their macroscopic motion and mechanical response[8]. Multibody dynamics analysis is a numerical computational technique primarily used to study the kinematic and dynamic behavior of electromechanical systems composed of multiple rigid or deformable bodies. In this type of analysis, the mechanical system is abstracted as a set of these bodies. The motion and force situations are then solved by modeling the interaction forces and motion states between them[9]. Fluid dynamics analysis mainly utilizes FEM to perform numerical simulations and calculations of physical phenomena involving gases or liquids. This analysis is used to predict and evaluate parameters such as fluid motion, pressure distribution, and temperature distribution under certain conditions, as well as to study the effects of fluids on the surrounding environment and objects[10]. In addition, FEM can also be used in medical fields, such as studying the influence of the geometric shape of antibiotic cement hip joint spacer(CS) on its mechanical properties[11], analyzing biomechanical stress distribution in bones and dental implants in dental implantology[12], and application of middle ear prosthesis design and personalization[13].

In recent decades, finite element analysis and finite element simulation have become increasingly prevalent in research in material science, including areas such as the stochastic homogenization problem (estimating the influence of microscopic random variation of heterogeneous material on a homogenized material property)[14], design of ultrasonic transducers[15], application in smart materials and structures[16], calculation of contour integral parameters for pipe welds[17], and design of multi-functional structures made of honeycomb sandwich materials. Although the finite element method (FEM) has been widely applied in various fields and can provide intuitive force distribution diagrams for research objects, there is still no objective criterion for evaluating the results of FEM analysis. For example, after making minor modifications to a component, deformation analysis was conducted using FEM. However, the deformation distribution map provided did not show significant differences to the naked eye. So, how can we determine whether the modified component has improved or deteriorated in terms of deformation distribution?

In this paper, we give a kind of quasi-distribution evaluation criteria based on piecewise Bézier curves to address the issue of the inability to objectively evaluate finite element models. During the optimization design of mechanical parts, finite element modeling is performed on their stress deformation, and the mesh node shape variable values are converted into distribution histogram data for piecewise Bézier curve fitting. Being dealt with area normalization method, the fitting curve could be regarded as a kind of probability density function (PDF), whose variance could be used to evaluate the finite element modeling results, and the situation with the minimum variance is the optimal choice for overall deformation.

In CAGD (computer aided geometric design) the Bernstein-Bézier form is the usual form of representing a polynomial curve[18]. Bézier curves have optimal shape preserving properties, and a Bézier curve of order $n$ is evaluated by the de Casteljau algorithm with a computational cost of $O(n^2)$ elementary operations[19]. But Bézier curve has a shortcoming, it is a kind of integral curve[20], thus means changing of even one control point, all the points in the curve will be changed. So, in this paper we will use piecewise Bézier curve to do the simulation work.

## 2  Theoretical considerations

When designing parts for industrial products, the first step is to create a 3D model of the parts and make several minor modifications according to actual needs. Then, based on the results of finite element analysis, the optimal design is selected. Firstly, we established a part model using the UGnx software system, which is made of 45 Mn steel with a tensile strength of 620 MPa, yield strength of 375 MPa, elongation of 15%, section shrinkage of 40%, and hardness of 240 HBW (see table 1). The shape of the part is a rectangular prism measuring 100 mm x 50 mm x 30 mm, with holes drilled on a square surface and the punched surface set to the front. Then, the coordinate system (unit: mm) is constructed with the lower left vertex of the selected square as the origin, and the sides connected to the vertex align with the positive directions of the x-axis and y-axis, respectively. The point (25, 75) is designated as point A, which is the drilling point (see Figure 1). We assume there are multiple options for the location of the stamping point, with point A as the center. The other stamping points are set along five directions: (-1,1), (-1,0), (-1,-1), (0,-1), and (1,-1). There are three selection positions in each direction, with the center distribution of the stamping point gradually moving away from point A, and the distance increasing by 1 mm (see Figure 2). The hole diameters of all perforated points are 10 mm, and the hole depth runs through the entire steel.

**Tab.1.** Detailed parameter values of the model

| MODEL PARAMETER | DATA |
|---|---|
| MODEL MATERIAL | 45Mn Steel |
| TENSILE STRENGTH | 620MPa |
| YIELD STRENGTH | 375MPa |
| ELONGATION | 15% |
| SECTION SHRINKAGE | 40% |
| HARDNESS | 240HBW |

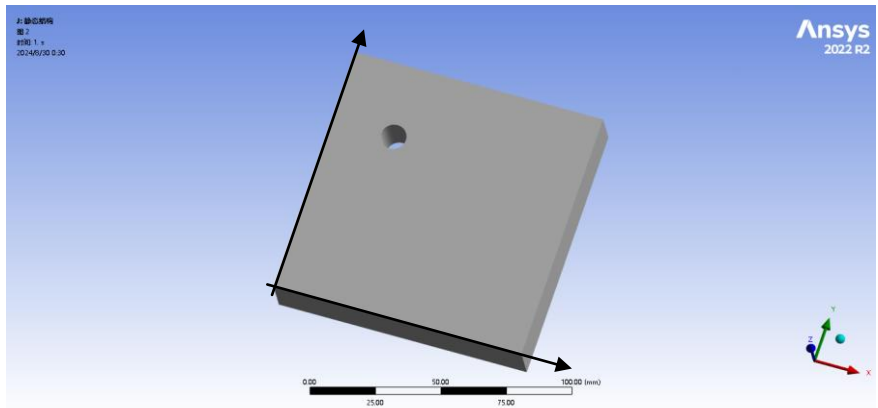

**Fig. 1.** 3-D graph of the part model

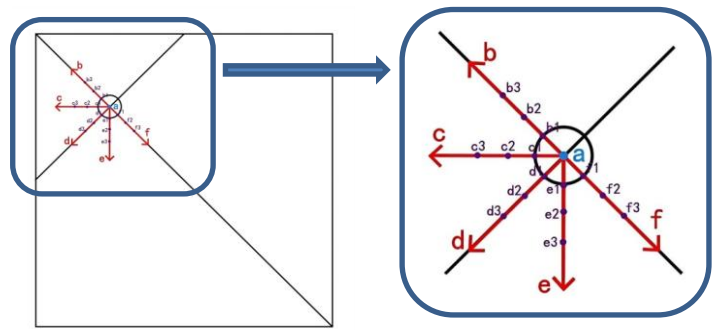

**Fig. 2.** Selection of the central position of drilling point

To conduct stress deformation analysis on the parts mentioned earlier, the part model is imported into ANSYS 2022 R2 for finite element analysis. To ensure the reliability of the model and the uniform distribution of the load on the stress surface during this static analysis process, the static analysis type is selected. First, we create a one-to-one corresponding static analysis system for all models and import the models into the geometric structure. Next, a linear mesh grid is created to make the distribution of nodes in the model grid more uniform. The grid accuracy is set to 4.0 mm, resulting in approximately 4950 nodes per model (see Figure 3). With this accuracy, the reliability of the solution exceeds 98%. In ANSYS, model analysis depends on the creation of grid nodes, and all applied loads and boundary conditions are applied to the nodes. The displacement of nodes corresponds to the deformation of the model.

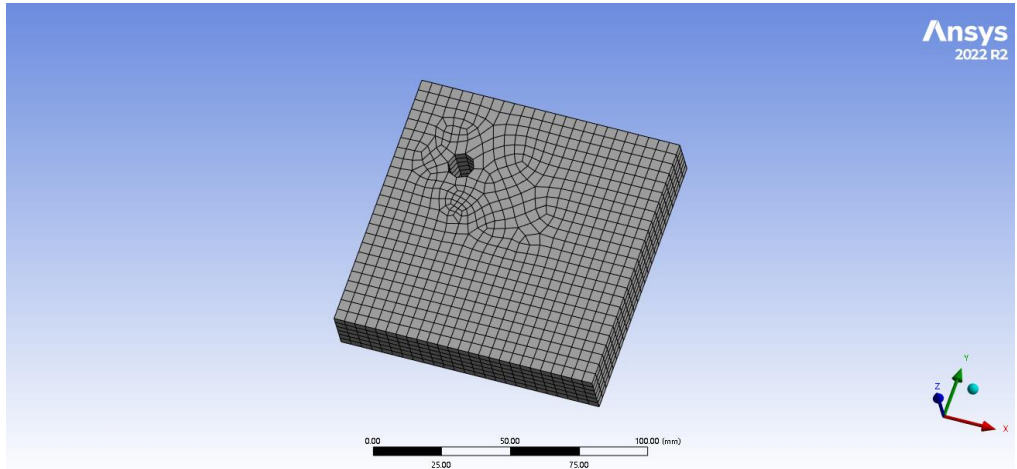

**Fig. 3.** Linear mechanical grid of model

In order to prevent the force applied to the model from being transformed into kinetic energy of the model without deformation, a fixed support is applied to the bottom of the model to prevent it from moving in any direction, (see figure 4).

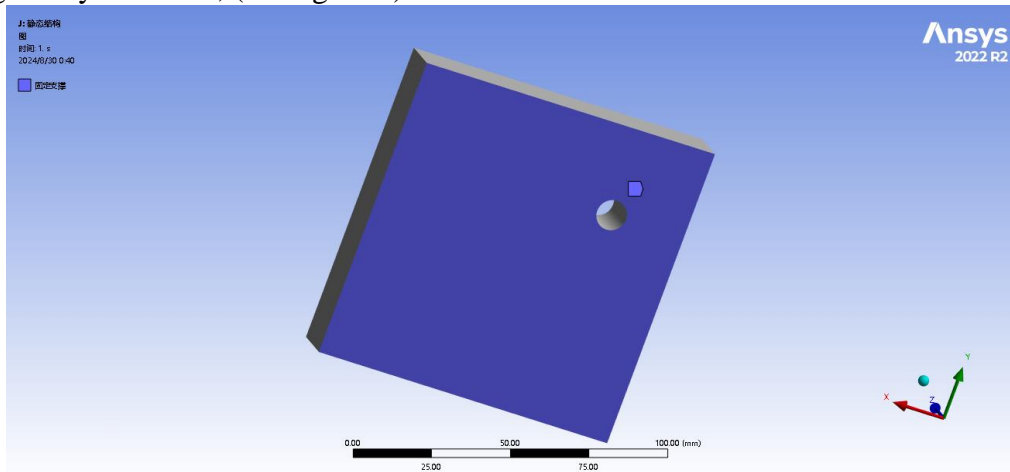

**Fig. 4.** Bottom fixing of model

After the position of the model is fixed, an evenly distributed vector force is applied to the front of the model, whose direction points vertically to the bottom and the magnitude of the load is 10000N, (see figure 5).

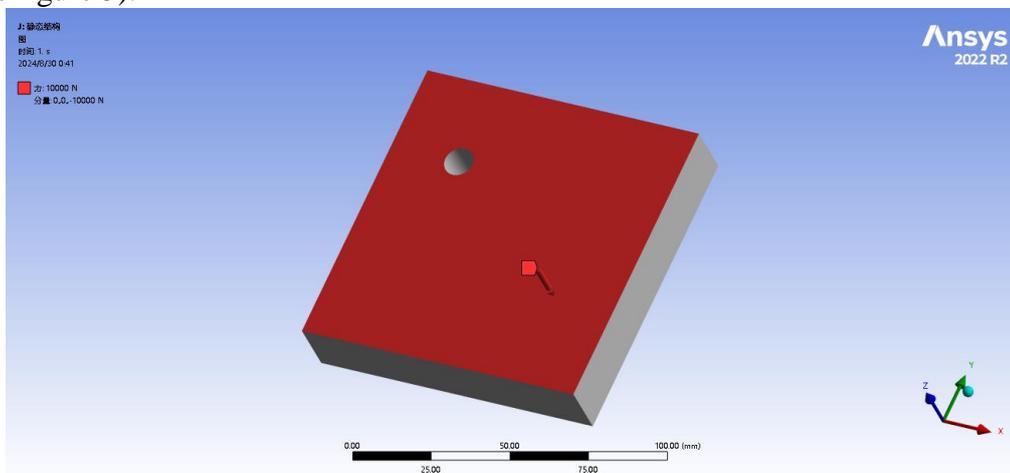

**Fig. 5.** Vector force applied to model

Finally, analyze the deformation of the model under the compression of the load and carry out finite element analysis on the model to make the stress distribution map of the model under stress, of which the different colors represent different stress levels, (see figure 6).

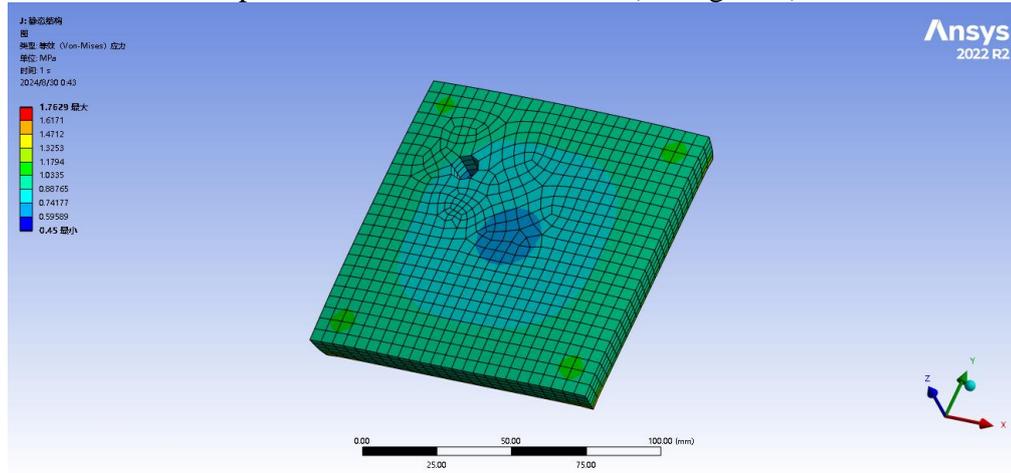

**Fig. 6.** Result of finite element analysis on the model

Based on the previous assumptions, there are 16 cases of finite element analysis results totally (demonstrated in figure 2), which are denoted as: $D_a$, $D_{b\_1}$, $D_{b\_2}$, $D_{b\_3}$, $D_{c\_1}$, $D_{c\_2}$, $D_{c\_3}$, $D_{d\_1}$, $D_{d\_2}$, $D_{d\_3}$, $D_{e\_1}$, $D_{e\_2}$, $D_{e\_3}$, $D_{f\_1}$, $D_{f\_2}$, $D_{f\_3}$. The deformation distribution map provided by those 16 scenarios could not show significant differences, we will evaluate them using the method provided below and select the optimal result.

## 2.1 Histogram distribution

Denote the total number of nodes in a single model as $N_T$, then the number of nodes for all finite element analysis results are shown in table 2.

**Tab.2.** The number of nodes for all finite element analysis results

| CASE | $N_T$ | CASE | $N_T$ | CASE | $N_T$ | CASE | $N_T$ |
|---|---|---|---|---|---|---|---|
| $D_a$ | 5068 | $D_{c\_1}$ | 5117 | $D_{d\_2}$ | 4935 | $D_{e\_3}$ | 4963 |
| $D_{b\_1}$ | 5106 | $D_{c\_2}$ | 5068 | $D_{d\_3}$ | 4991 | $D_{f\_1}$ | 5215 |
| $D_{b\_2}$ | 5180 | $D_{c\_3}$ | 4963 | $D_{e\_1}$ | 5075 | $D_{f\_2}$ | 5054 |
| $D_{b\_3}$ | 5047 | $D_{d\_1}$ | 5033 | $D_{e\_2}$ | 4977 | $D_{f\_3}$ | 5159 |

When we build the finite element analysis on the model to make the stress distribution map of the model under stress, the program assigns numbers to each node on the grid in turn, calculates the displacement of each node in $X$, $Y$ and $Z$ directions. Let $V_i = \begin{pmatrix} V_i^X & V_i^Y & V_i^Z \end{pmatrix}$ denote displacement value of node $i$, obviously $i = 1, 2, \cdots N_T$. For ease of calculation, use 1-norm to measure the value of vector $V_i$, thus $\|V_i\|_1 = |V_i^X| + |V_i^Y| + |V_i^Z|$, $i = 1, 2, \cdots N_T$.

Then we build histogram distribution of given finite element analysis model as following:
1) Assume $V_{\min} = \min\limits_{1 \leq i \leq N_T} \|V_i\|_1$, $V_{\max} = \max\limits_{1 \leq i \leq N_T} \|V_i\|_1$, $\xi = (V_{\max} - V_{\min})/350$, divided the interval $(V_{\min}, V_{\max})$ into 350 equal parts.

Assume $P_k(\|V_i\|_1) = \begin{cases} 1, & \text{if } \|V_i\|_1 \in (V_{\min} + (k-1)\xi, F_{\min} + k\xi) \\ 0, & \text{else} \end{cases}$; $k = 1, 2, \cdots 350$; $i = 1, 2, \cdots N_T$.

2) Let $\bar{P}_k = \dfrac{\sum_{i=1}^{N_T} P_k(\|V_i\|_1)}{N_T}$ $k = 1, 2, \cdots 350$, then $\sum_{k=1}^{350} \bar{P}_k = 1$. That means $\bar{P}_k$ could be regard as a kind of probability distribution. We can call that histogram distribution, which could be considered as a dimensionless signal with the length of 350.

## 2.2 Sub-section Bézier curve

For those histogram distribution data $\bar{P}_k$, it could be simulated with a function $P(x)$ which has the properties of probability density function (PDF). Thus the histogram distribution data $\bar{P}_k$ could be analysis at the angle of probability. In this paper, we simulate $\bar{P}_k$ with piecewise Bézier curve, which is more flexible in curve modeling. piecewise Bézier curve is a kind of parameter curve (the parameter denote as $t \in [0,1]$), in 5 degree situation, it has 11 bases, denote as $N_i$, $i = 0, 1, \cdots 10$:

$$N_i(t) = \begin{cases} \binom{5}{i}\left(1-\dfrac{t}{\omega}\right)^{5-i} \cdot \left(\dfrac{t}{\omega}\right)^i & t \in [0, \omega) \\ 0 & t \in [\omega, 1] \end{cases} \quad i = 0,1,2,3,4, \quad N_5(t) = \begin{cases} \left(\dfrac{t}{\omega}\right)^5 & t \in [0, \omega) \\ \left(\dfrac{1-t}{1-\omega}\right)^5 & t \in [\omega, 1] \end{cases},$$

$$N_i(t) = \begin{cases} 0 & t \in [0, \omega) \\ \binom{5}{i-5}\dfrac{(1-t)^{10-i}(t-\omega)^{i-5}}{(1-\omega)^5} & t \in [\omega, 1] \end{cases} \quad i = 6,7,8,9,10,$$ where $\omega \in (0,1)$ is the segmentation point. Based on the models provided in this paper, it was found through calculation that $\omega = 0.608$ is the optimal choice.

Assume $C_i, i = 0, 1, \cdots 10$ are the points in two-dimensional plane, then the definition of 5 degree piecewise Bézier curve is:

$$B(t) = \sum_{i=0}^{10} N_i(t) C_i \quad t \in [0,1]. \tag{1}$$

Where $C_i$ are called control points of the piecewise Bézier curve defined by equation (1).

## 3 Simulation result

Simulation means, for given data $\bar{P}_k$, find a piecewise Bézier curve $B(t) = \sum_{i=0}^{10} N_i(t) C_i$ to simulate those data. The most important part is to calculate unknown control points $C_i$.

### 3.1 Least square approximation method

In this paper, the least square approximation is used to deal with this problem. First, parameterize those data $\bar{P}_k$ by cumulative chord length parameterization method to match a parameter $t_k$ for every $\bar{P}_k$, thus we get a parametric sequence $0 = t_1 < \cdots t_k < \cdots < t_{1000} = 1$. Second, build a vector equation group which has 350 equations to solve the unknown control points $C_i$:

$$B(t_k) = \sum_{i=0}^{10} N_i(t_k) C_i = \bar{P}_k, \quad k = 1, 2, \cdots 350. \tag{2}$$

Equation (2) could be solved as:

$$\phi^T\phi\begin{bmatrix}C_0\\C_1\\\vdots\\C_{10}\end{bmatrix}=\phi^T\begin{bmatrix}\bar{P}_0\\\bar{P}_1\\\vdots\\\bar{P}_{350}\end{bmatrix}, \text{ where }\phi=\begin{bmatrix}N_0(t_1)&N_1(t_1)&\cdots&N_{10}(t_1)\\N_0(t_2)&N_1(t_2)&\cdots&N_{10}(t_2)\\\vdots&\vdots&&\vdots\\N_0(t_{350})&N_1(t_{350})&\cdots&N_{10}(t_{350})\end{bmatrix}. \quad (3)$$

In order to decide the segmentation point $\omega$, we need to figure out how to evaluate the goodness of a simulation result. As a parameter simulation curve, $B(t)$ need to be discretized into the standard discrete signal $\bar{B}_k$ to match with $\bar{P}_k$. Divide parameter interval $(0,1)$ into $n$ uniform parts ($n \gg 1000$), then get a discrete signal sequence $B_k=(B_k^x,B_k^y), k=1,2,\cdots n$. The x-coordinate of $\bar{P}_k$ is $k$, then we get $\bar{B}_k = B_{k'}^y$, if $B_{k'}^x = \max_i\{B_i^x < k\}$. Then we evaluate the goodness of the simulation result with mean square deviation (MSE), calculated as:

$$MSE(\bar{B}_k,\bar{P}_k)=\frac{\sum_{i=1}^{350}(\bar{B}_i-\bar{P}_i)^2}{350}. \quad (4)$$

Denote the standard discrete signal $\bar{B}_k$ got under the segmentation point $\omega'$ as $\bar{B}_k^{\omega'}$, then the best segmentation point $\omega$ could be decided like:

$$\omega=\min_{\omega'\in(0,1)} MSE(\bar{B}_k^{\omega'},\bar{P}_k) \quad (5)$$

### 3.2 Quasi-distribution

The simulation curve $B(t)$ is an approximation to the histogram distribution $\bar{P}_k$, thus the integral of $B(t)$ should equal 1. Apparently, that is not always come to existence. But we can fulfill it through an adjustment factor $\gamma$, which is reciprocal of the area enclosed by the simulation curve $B(t)$ and the x-axis. Let $\bar{B}(t)=\gamma B(t)$, correspondingly, $\bar{B}(t)$ satisfy the property of probability density function (PDF). But its expression is not like any existed PDF, then, we can call it quasi-distribution. Since all the existed probability density functions are based on the independence of random variable sequence which is not always come into existence, so the quasi-distribution mentioned before is meaningful.

### 3.3 Experimental results

To the finite element analysis results of 16 cases: $D_a$, $D_{b\_1}$, $D_{b\_2}$, $D_{b\_3}$, $D_{c\_1}$, $D_{c\_2}$, $D_{c\_3}$, $D_{d\_1}$, $D_{d\_2}$, $D_{d\_3}$, $D_{e\_1}$, $D_{e\_2}$, $D_{e\_3}$, $D_{f\_1}$, $D_{f\_2}$, $D_{f\_3}$ (as in figure2), we get the quasi-distribution result show as figure3:

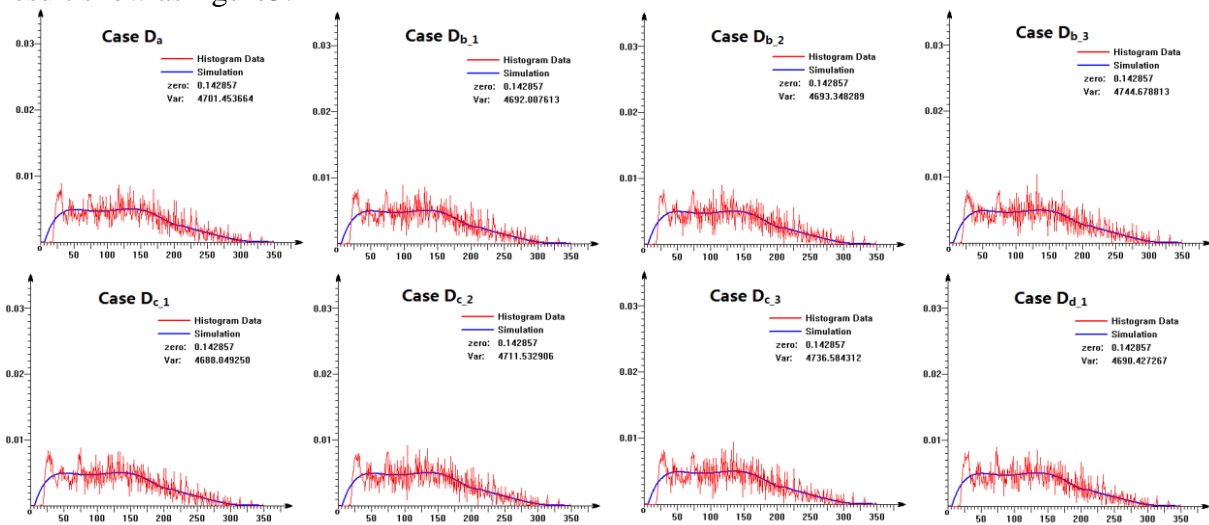

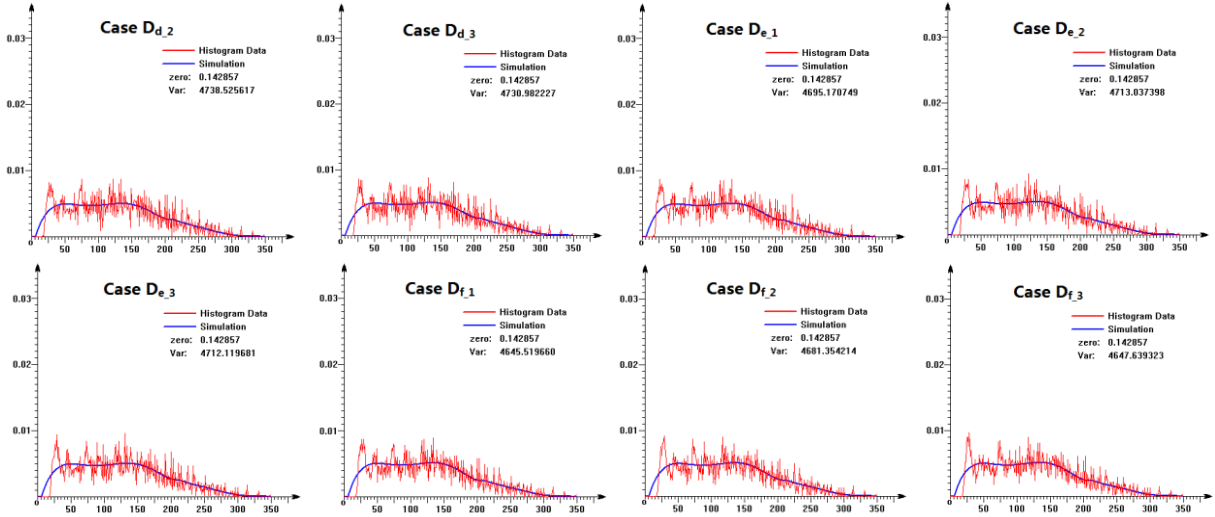

**Fig. 7.** Histogram data and quasi-distribution simulation

In each scenario of figure7, the red curve represents the value of the histogram data, and the blue curve represents the corresponding quasi-distribution simulations curve. And "zero" represents the proportion of nodes with zero deformation to the total number of nodes $N_T$. By calculating the variance of the quasi-distribution simulations curve (actually, the variance of signal $\bar{B}_k$ got from the quasi-distribution simulations curve), we can tell the case $D_{f\_1}$ could be regarded as the best case, because of the minimum variance, which means the optimal choice for overall deformation. And the case $D_{f\_3}$ is also a good choice, with a variance value only slightly larger than case $D_{f\_1}$. Note that both options are located in the same direction (see figure 2), indicating the effectiveness of the method presented in this paper.

## 4    Conclusion

In this paper, we present a new method for the optimization problem of component design, called quasi-distribution appraisal based on piecewise Bézier curves, to evaluate the finite element analysis of deformation for a given part. By simulating the distribution histogram data generated from the numerical deformation values of finite element analysis derived from different design cases with piecewise Bézier curves and area normalization methods, the simulation curve can be interpreted as a probability density function, termed quasi-distribution. The results of numerical simulation reveal that the histogram distribution data exhibit a distinct intrinsic pattern, which can be demonstrated by the quasi-distribution proposed in this paper. Subsequently, the corresponding variance can be used to assess the finite element analysis results and determine the best design case.

**Acknowledgement**

We would like to dedicate this manuscript to Pro. Ping Xi, whose commitment to this project and to scientific discovery were crucial in propelling this study forward. Dr. Jixing Li will be remembered for his intelligence, drive, and love for science and for all of the lives that he touched during his impactful career.

**References**

1. Amaze, C., Kuharat, S., Bég, O.A. et al, Finite element stress analysis and topological optimization of a commercial aircraft seat structure[J]. European Mechanical Science, 8(2), pp. 1–17, 2024.
2. Farkas, L., Moens, D. & Vandepitte, D., The interval finite element method for static structural analysis[C]. Proceedings of the Third M.I.T Conference on Computational Fluid and Solid Mechanics, pp. 202–205, 2005.
3. Yu, S.B. & Wang, H., Investigation of circumferential mode frequencies of circular cylindrical

shells of stator in electric motor[J]. Electric Machines and Control, 18(6), pp. 102–107, 2014.
4. Kim, S.W., Song J.W., Hong J.P. & Kim, H.J. et al, Finite element analysis and validation of wind turbine bearings[J]. Energies, 17(2), pp. 692, 2024.
5. Maier, W., Moehring, H.C. & Wunderle, R., Augmented reality to visualize a finite element analysis for assessing clamping concepts[J]. the International Journal of Advanced Manufacturing Technology, 133(5-6), pp. 1–10, 2024.
6. Wang, Q., Li, Z.F., Wang, F. & Li H., Fatigue life prediction of power transmission towers based on nonlinear finite element analysis[J]. Applied Mathematics and Nonlinear Sciences, 9(1), pp. 1–16, 2023.
7. Wang, L., Yue, C.X., Liu, X.L. et al, Conventional and micro scale finite element modeling for metal cutting process: a review[J]. Chinese Journal of Aeronautics, 37(2), pp. 199–232, 2024.
8. Müzel, S.D., Bonhin, E.P., Guimarães, N.M. & Guidi, E.S., Application of the finite element method in the analysis of composite materials: a review[J]. Polymers, 12(4), pp. 818, 2020.
9. Shao, H.Y., Li, D.C., Kan, Z., Xiang, J.W. & Wang, C.S., Analysis of catapult-assisted takeoff of carrier-based aircraft based on finite element method and multibody dynamics coupling method[J]. Aerospace, 10(12), pp. 1005, 2023.
10. Huang, Y., Xiao, Q., Idarraga, C. et al, Novel computational fluid dynamics-finite element analysis solution for the study of flexible material wave energy converters[J]. Physics of Fluids, 35(8), pp. 083611, 2023.
11. Warinsiriruk, E., Thongchuea, N., Pengrung, N. et al, Inverstigation of mechanical behavior on the cement hip spacer geometry under finite element method and compression load test[J]. Journal of Orthopaedics, 47(11), pp. 115–121, 2024.
12. Qiu, P.P., Cao, R.K., Li, Z.Y. & Fan, Z., A comprehensive biomechanical evaluation of length and diameter of dental implants using finite element analyses: a systematic review[J]. Heliyon, 10(6), pp. 26876, 2024.
13. Mohseni-Dargah, M., Pastras, C., Mukherjee, P., Khajeh, K., & Asadnia, M., Enhancing ossicular chain reconstruction through finite element analysis and advanced additive manufacturing: a review[J]. Bioprinting, 38(1), pp. e00328, 2024.
14. Fosness, E., Guerrero, J., Qassim, K. & Denoyer, S.J., Recent advances in multi-functional structures. 2000 IEEE Aerospace Conference Proceedings, IEEE, pp. 23–28, 2000.
15. Kim, J., Kim, H.S, Finite element analysis of piezoelectric underwater transducers for acoustic characteristics[J]. Journal of Mechanical Science and Technology, 23(2), pp.452-460, 2010.
16. Sakata, S.I. & Torigoe I., A successive perturbation-based multi-scale stochastic analysis method for composite materials[J]. Finite Elements in Analysis and Design, 102-103(1), pp. 74–84, 2015.
17. Benjeddou, A., Advances in piezoelectric finite element modeling of adaptive structural element: a survey[J]. Computer and Structures, 76(2), pp. 347–363, 2000.
18. Wang, Y.D., Hu, B.F., Xi, P. & Xue, W., A note on variable upper limit integral of Bézier curve[J]. Advanced Science Letters, 4(8-10), pp. 2986–2990, 2011.
19. Zhao, Q.L., Chen, J.J., Feng, X.B. & Wang, Y.D., A novel Bézier LSTM model: a case study in corn analysis[J]. Mathematics, 12(15), pp. 2308, 2024.
20. Zhao, Q.L., Lu, Z.H. & Wang, Y.D., Analyzing the data of COVID-19 with quasi-distribution fitting based on piecewise B-spline curves[J]. COVID, 2(2), pp. 175-196, 2022.